\documentclass{amsart}
\usepackage{hyperref}
\usepackage{amssymb} 
\usepackage{comment, enumerate}
\usepackage{comment}
\usepackage[all,cmtip]{xy}
\usepackage{url}

\DeclareMathOperator{\Sp}{Sp}
\DeclareMathOperator{\GSp}{GSp}
\DeclareMathOperator{\Stab}{Stab}

\DeclareMathOperator{\ord}{ord}

\DeclareMathOperator{\GL}{GL}

\DeclareMathOperator{\Gal}{Gal}

\newcommand{\Q}{{\mathbb Q}}
\newcommand{\Z}{{\mathbb Z}}
\newcommand{\F}{{\mathbb F}}

\def\C{{\ensuremath{\mathbb{C}}}}

\def\O{{\ensuremath{\mathcal{O} }}}

\def\F{{\ensuremath{\mathbb{F}}}}
\def\Flb{{\ensuremath{\overline{\mathbb{F}}_\ell}}}

\def\p{{\ensuremath{\mathfrak{p}}}}
\def\s{{\ensuremath{\mathfrak{s}}}}

\def\triv{{\ensuremath{\mathbf{1}}}}

\DeclareMathOperator{\Frob}{Frob}

\DeclareMathOperator{\Id}{Id}

\DeclareMathOperator{\Ind}{Ind}

\DeclareMathOperator{\Jac}{Jac}

\begin{document}

\newtheorem{lem}{Lemma}[section]
\newtheorem{thm}[lem]{Theorem}
\newtheorem{prop}[lem]{Proposition}
\newtheorem{alg}[lem]{Algorithm}
\newtheorem{cor}[lem]{Corollary}
\newtheorem{conj}{Conjecture}
\newtheorem*{conj*}{Conjecture}

\theoremstyle{definition}
\newtheorem{dfn}[lem]{Definition}
\newtheorem{rmk}[lem]{Remark}
\newtheorem*{ex}{Example}

\theoremstyle{remark}

\newtheorem{ack}{Acknowledgement}

\title[Constructing hyperelliptic curves]{Constructing hyperelliptic curves with surjective Galois representations}
\author{Samuele Anni and Vladimir Dokchitser}

\address{Institut de Math\'ematiques de Marseille, Universit\'e d'Aix-Marseille, 163, avenue de Luminy, Case 907, 13288 Marseille Cedex 9, France}
\email{samuele.anni@gmail.com}

\address{Department of Mathematics, King's College London, Strand, London, WC2R 2LS, United Kingdom}
\email{vladimir.dokchitser@kcl.ac.uk}

\keywords{Galois representations, abelian varieties, hyperelliptic curves, inverse Galois problem, Goldbach's conjecture}
\subjclass[2010]{Primary 11F80, Secondary 12F12, 11G10, 11G30}
\begin{abstract}
In this paper we show how to explicitly write down equations of hyperelliptic curves over $\Q$ such that for all odd primes $\ell$ the image of the mod~$\ell$ Galois representation is the general symplectic group. 
The proof relies on understanding the action of inertia groups on the $\ell$-torsion of the Jacobian, including at primes where the Jacobian has non-semistable reduction.
We also give a framework for systematically dealing with primitivity of symplectic mod~$\ell$ Galois representations. 

The main result of the paper is the following. Suppose $n=2g+2$ is an even integer that can be written as a sum of two primes in two different ways, with none of the primes being the largest primes less than $n$ (this hypothesis is expected to hold for all $g\neq 0,1,2,3,4,5,7$ and $13$).
Then there is an explicit $N\in\Z$ and an explicit monic polynomial $f_0(x)\in \Z[x]$ of degree $n$, such that the Jacobian $J$ of every curve of the form $y^2=f(x)$ has 
$\Gal(\Q(J[\ell])/\Q)\cong \GSp_{2g}(\F_\ell)$ for all odd primes $\ell$ and $\Gal(\Q(J[2])/\Q)\cong S_{2g+2}$, whenever $f(x)\in\Z[x]$ is monic with $f(x)\equiv f_0(x) \bmod{N}$ 
and with no roots of multiplicity greater than $2$ in $\overline{\F}_p$ for any $p\nmid N$.
\end{abstract}
\maketitle
\vspace{-1cm}
\tableofcontents
\vspace{-1.2cm}
\section{Introduction}

A classical method for showing that the group $\GL_2(\F_\ell)$ is a Galois group over $\Q$ is by realising it as the Galois group of the field generated by the $\ell$-torsion on an elliptic curve. 
One can similarly try to construct the general symplectic group $\GSp_{2g}(\F_\ell)$ as the Galois group associated to the $\ell$-torsion of a $g$-dimensional abelian variety.
The main difficulty is that it is much harder to write down explicit abelian varieties and then verify that the Galois group obtained is not a proper subgroup of $\GSp_{2g}(\F_\ell)$. 
This approach has been successfully used in a number of works, including \cite{zarhin79}, \cite{zarhin_silv_05}, \cite{hall08}, \cite{zarhin10}, \cite{hall11}, \cite{reyna_vila_tame}, \cite{reyna_kappen},  \cite{arias15}, 
which realise $\GSp_{2g}(\F_\ell)$ for every (odd) prime $\ell$ using Jacobians of hyperelliptic curves, and show that one curve often realises $\GSp_{2g}(\F_\ell)$ for all sufficiently large $\ell$. 
More recently \cite{LandesmanSwaminathanTaoXu} gave a non-constructive proof that many hyperelliptic curves realise $\GSp_{2g}(\F_\ell)$ for all odd primes $\ell$, and 
\cite{Dieulefait4}, \cite{Zywina15}, \cite{annilemossiksek} who exhibited explicit curves of genus $2$ and $3$ with this property. 
There has also been numerical work \cite{reyna_armana_karemaker_rebolledo_thomas_vila}, \cite{annilemossiksek} investigating the Galois images of Jacobians of curves, and work on general abelian varieties \cite{Lombardo}.

The main contribution of the present article to this topic is an explicit construction of hyperelliptic curves, such that for every prime $\ell$ their Jacobian has maximal mod $\ell$ Galois image; 
in other words hyperelliptic curves whose Jacobian $J$ has $\Gal(\Q(J[\ell])/\Q)\cong\GSp_{2g}(\F_\ell)$ for every odd prime $\ell$, and isomorphic to the symmetric group $S_{2g+2}$ for $\ell=2$. 

The key new tool is a way of controlling the action of inertia groups on $J[\ell]$ at places where $J$ has non-semistable reduction. 
Our approach does, however, require a rather unorthodox constraint on the dimension $g$: we need the even integer $2g+2$ to satisfy Goldbach's conjecture! 
In fact, we require it to satisfy the following somewhat stronger statement, which appears to hold\footnote{We have made a quick numerical check: the property holds for all other $g<10^7$.} for every genus $g\neq 0,1,2,3,4,5,7,13$:

\begin{conj*}[Double Goldbach]
Every positive even integer $n$ can be written as a sum of two primes in two different ways with none of the primes being the largest prime less than $n$, except for $n= 0,2,4,6,8,10,12, 16,28$.
\end{conj*}

The result on Galois images of hyperelliptic curves that we obtain is the following:

\begin{thm} \label{final:congruence}
Let $g$ be a positive integer such that $2g+2=q_1+q_2=q_4+q_5$ for some primes $q_i$, with $\{q_1,q_2\}\neq \{q_4, q_5\}$, and that there is a further prime $q_3$ with $q_1,q_2,q_4,q_5<q_3<2g+2$. Then there exist an explicit $N\in \Z$ and an explicit $f_0\in \Z[x]$ monic of degree $2g+2$ such that if
\begin{enumerate}
 \item $f(x)\equiv f_0\bmod{N}$, and
 \item $f(x) \bmod{p}$ has no roots of multiplicity greater than $2$ for all primes $p\nmid N$,
\end{enumerate}
then $\Gal(\Q(J[\ell])/\Q)\cong \begin{cases} 
\GSp_{2g}(\F_\ell) \mbox{ for all primes }\ell\neq 2,\\
S_{2g+2} \mbox{ for }\ell=2,
\end{cases}$ \hspace{-4.5mm} where ${J{=}\Jac(y^2{=}f(x))}$.
\end{thm}

\noindent 
See Theorem~\ref{final:congruence_th} for the explicit description of $N$ and $f_0(x)$, and Remark~\ref{explicit_notriple} for an explanation on how to find explicit curves satisfying hypothesis $(2)$. 
An explicit example for $g=6$ is given in Section~\ref{section_example}.

For the exceptional genera $g=1,2,3,4,5,7,13$ our method still makes it possible to construct hyperelliptic curves with maximal image at all but a small number of primes, e.g. all primes except $\ell=5, 11, 13$ when $g=7$ (see Remark~\ref{exceptional_list}). 

It is worth noting that hypotheses $(1)$ and $(2)$ in the above theorem are satisfied by a positive density of monic polynomials $f(x)$, for example see \cite[Theorem~$1.5$]{Bhargava}. 
In particular, this shows that the hyperelliptic curves of genus $g$ with maximal mod $\ell$ Galois image for every prime $\ell$ have a positive (lower) density among all hyperelliptic curves of genus $g$.

\bigskip

Throughout the paper we work with hyperelliptic curves 
$$
C: y^2 = f(x)
$$
and write
$$
 J =\Jac(C)
$$
for their Jacobian. The layout is as follows.

In Section~\ref{section_inertia} we examine the Galois representations $H_{\acute{e}t}^1(C,\Q_\ell)$ and $J[\ell]$ as representations of local Galois groups. 
For the ``$\ell\neq p$'' theory we use the method of clusters, recently introduced in \cite{DDMM}. 
For ``$\ell=p$'' we restrict our attention to primes $p$ of semistable reduction, and use the description given by the theory of fundamental characters. 
We also give a simple criterion guaranteeing that a local Galois group contains a transvection in its action on $J[\ell]$.

In Section~\ref{section_irred} we develop criteria for the representations $H_{\acute{e}t}^1(C,\Q_\ell)$ and $J[\ell]$ to be globally irreducible. 
This is the place where the ``double Goldbach'' hypothesis enters. The reason for it is that we cannot always guarantee that the above representations are locally irreducible.
In fact, this appears to be a genuine obstruction: for example, there is no $17$-dimensional abelian variety $A/\Q_p$ and primes $p>35$ and $\ell\neq p$, such that $T_\ell(A)\otimes\overline{\Q}_\ell$ is an irreducible 
$\Gal(\overline{\Q}_p/\Q_p)$-module\footnote{The hypotheses ensure that $A$ has potentially good reduction and the inertia group at $p$ acts tamely and semisimply on $T_\ell(A)$, through a cyclic quotient of order $k$, say. Any element of the inertia group must have rational trace on $T_\ell(A)$ (as the trace is independent of $\ell$), so irreducibility forces the eigenvalues of a generator of the image of inertia to be precisely the set of $k$-th roots of unity. Hence $34 = \dim T_\ell(A)\otimes\overline{\Q}_\ell = \varphi(k)$, which is impossible.}.
However, when $2g+2$ is a sum of two primes (other than $\ell$), we are able to force the local representation to have at most two irreducible constituents. To guarantee global irreducibility we then use conditions of this kind at several places. 
The reason that we require ``double Goldbach'' rather than the classical Goldbach's conjecture, is so that we can treat the case when $\ell$ is one of the Goldbach prime summands of $2g+2$ by using the other pair of primes.

In Section~\ref{section_prim} we develop criteria for the representations $H_{\acute{e}t}^1(C,\Q_\ell)$ and $J[\ell]$ to be primitive. 
The basic method essentially follows that of Serre (see e.g.~\cite[Theorem~4]{Mazur78} or~\cite[\S 1]{annisiksek_unif}) or, more recently \cite{annilemossiksek}: 
we ensure that every inertia group acts trivially on every possible Galois stable partition of the representation and then invoke the Hermite-Minkowski theorem to deduce that no such partition exists. 
We formalise this approach by introducing ``quasi-unramified'' representations, and develop the
necessary conditions on $f(x)$ that make the argument work (``admissible'' and ``$p$-admissible'' polynomials). Unlike \cite{annilemossiksek}, it is important for us to allow for curves with non-semistable reduction.

In Section~\ref{section_surj} we recall the classification of subgroups of $\GSp_{2g}(\F_\ell)$ of \cite{hall08} and \cite{disawi} and rephrase it as a criterion for $J$ to have $\Gal(\Q(J[\ell])/\Q)\cong\GSp_{2g}(\F_\ell)$. 
As in many previous works, the basic rule is that if the action is irreducible, primitive and contains a transvection, then the Galois representation has maximal image
(see Theorem~\ref{classification}).

In Section~\ref{section_maximal} we tie everything together to give a list of (essentially local) constraints that guarantee that $\Gal(\Q(J[\ell])/\Q)\cong\GSp_{2g}(\F_\ell)$ for every odd prime $\ell$, and that 
$\Gal(\Q(J[2])/\Q)\cong S_{2g+2}$; see Theorems~\ref{out:7} and~\ref{alldone}. 

In Section~\ref{section_congruence} we give explicit congruence conditions on the coefficients of $f(x)$ that ensure that the above list of constraints is satisfied, and prove the precise version of Theorem~\ref{final:congruence} (see Theorem~\ref{final:congruence_th}).

We end in Section~\ref{section_example} by working through the conditions in Theorem~\ref{final:congruence_th} for $g=6$, and constructing a hyperelliptic curve satisfying all the hypotheses.\\

\noindent{\bf Acknowledgments.} We would like to thank Adam Morgan and Tim Dokchitser for several useful discussions related to this work.
The first author was supported by EPSRC Programme Grant \lq LMF: L-Functions and Modular Forms\rq\  EP/K034383/1 during his position at the University of Warwick, and by DFG Priority Program SPP 1489 and the Luxembourg FNR during his positions at IWR, Heidelberg and at the University of Luxembourg. The second author is supported by a Royal Society University Research Fellowship.
Both authors would also like to thank the Warwick Mathematics Institute where most of this research was carried.

\subsection{Notation,  $t$-Eisenstein polynomials and type $t-\{q_1,{\dots},q_k\}$}
$\;$\\

\noindent
{\underline{{\bf Local setting}}. For a finite extension  $F$ of $\Q_p$ we write:
\begin{itemize}
\item $\pi_F$ for a fixed uniformizer of $F$;
\item $\O_F$ for the ring of integers of $F$;
\item $e_F$ for the ramification degree of $F$;
\item $\overline{F}$ for a fixed algebraic closure of $F$;
\item $F^{\mathrm{nr}}$ for the maximal unramified extension of $F$;
\item $v$ for a valuation on $\overline{F}$ normalized such that $v(\pi_F)=1$;
\item $\F$ for the residue field of $F$;
\item $G_F$ for the absolute Galois group $\Gal(\overline{F}/F)$;
\item $I_{F}$ for the inertia subgroup of $G_F$;
\item $\overline{g}(x)$ and $\overline{\alpha}$ for the reduction modulo $\pi$ of every $g(x)\in \O_F[x]$ and $\alpha\in \O_F$.
\end{itemize}

\begin{dfn}[$t$-Eisenstein polynomials]
\label{t-eis}
Let $t\geq 1$ be an integer. We say that a polynomial with $\O_F$-coefficients 
$$f(x)=x^m+a_{m-1}x^{m-1}+\dots+a_0$$ is $t$-\emph{Eisenstein} if $v(a_i)\geq t$ for all $i>0$ and $v(a_0)=t$. 
\end{dfn}

\begin{dfn}[Polynomials of type $t-\{q_1,{\dots},q_k\}$]
\label{type}
Let $q_1,{\dots},q_k$ be prime numbers and let $t\geq 1$ be an integer. Let $f(x)\in \O_F[x]$ be a monic squarefree polynomial. We say that $f(x)$ is of \emph{type} $t-\{q_1,{\dots},q_k\}$ if it can be factored as  
$$f(x)=h(x)\prod_{i=1}^k g_i(x-\alpha_i)$$ over $\O_F[x]$, for some $\alpha_i \in \O_F$ with $\overline{\alpha}_i \neq \overline{\alpha}_j$ for all $i\neq j$, 
where $g_i(x)$ is a $t$-Eisenstein polynomial of degree $q_i$ and $\overline{h}(x)$ is separable with $\overline{h(\alpha_i)}\neq 0$ for all $i$.

In other words, the monic polynomial $f(x)$ is a product of shifted $t$-Eisenstein polynomials of degrees $q_i, \dots q_k$ and linear polynomials, 
such that these polynomials have no common roots in the residue field. See Section~\ref{section_example} for explicit examples.
\end{dfn}

\noindent
{\underline{{\bf Global setting}}. For a number field $K$ we write:
\begin{itemize}
\item $\O_K$ for the ring of integers of $K$;
\item $\F_\p$ for the residue field of $K$ at a prime $\p$ of $K$;
\item $\overline{K}$ for a fixed algebraic closure of $K$;
\item $G_K$ for the absolute Galois group $\Gal(\overline{K}/K)$;
\item $I_{\p}=I_{K_\p}$ for the inertia subgroup at $\p$.
\end{itemize}

\begin{dfn}
Let $q_1,{\dots},q_k$ be prime numbers. Let $f(x)\in \O_K[x]$ be a monic squarefree polynomial. We say that $f$ is of \emph{type} $t-\{q_1,{\dots},q_k\}$ at $\p$ if 
$f(x)\in \O_{K_\p}[x]$ is of type $t-\{q_1,{\dots},q_k\}$. 
\end{dfn}

\noindent
{\underline{{\bf Roots of unity}}. Let $q$ be a positive integer, we will denote by $\zeta_{q}$ a primitive $q$-th root of unity. 
Throughout this article we will choose primitive roots of unity to form compatible systems, i.e.\ if $\zeta_q$ is a primitive $q$-th root of unity and $q=q'q''$ then $\zeta_q^{q''}=\zeta_{q'}$. 
In characteristic $\ell$ dividing $q$, we have $\zeta_{q}=\zeta_m$ where $\ell^r m=q$, with $(\ell, m)=1$.

\section{Inertia action on $J[\ell]$}\label{section_inertia}

The construction of hyperelliptic curves presented in this article will crucially rely on understanding the action of the inertia groups on the $\ell$-torsion of the Jacobian for every prime number $\ell$.
In this section, $F$ will be a local field of odd residue characteristic $p$. 
Let $$C:y^2=f(x)$$ be a genus $g$ hyperelliptic curve over $F$ with $f(x)\in \O_F[x]$ monic and squarefree, and let $J=\Jac(C)$. 
We will describe the inertia action on $J[\ell]$ in terms of $f(x)$.

\subsection{$J[\ell]$ when $\ell\neq p$: clusters}\label{section_inertia_clusters}
$\;$\\

\noindent
In this section we describe the inertia action on $J[\ell]$ when $\ell\neq p$. In particular, we will prove the following theorem:

\begin{thm}
\label{eigen:tame}
Suppose that $f(x)\in \O_F[x]$ has type $t-\{q_1,{\dots},q_k\}$ for odd primes $q_i\neq p$. 
Then for every $\ell\neq p$, the inertia group $I_F$ acts tamely on $H^1_{\acute{e}t}(C/\overline{F}, \Q_\ell)$ and on $J[\ell]$ through a quotient of order dividing $2\prod_i q_i$.
Moreover, the non-trivial eigenvalues (with multiplicity) of any generator $\tau$ of tame inertia are either
$$-\zeta_{q_1},-\zeta_{q_1}^2,\dots,-\zeta_{q_1}^{q_1-1},\dots,-\zeta_{q_k},-\zeta_{q_k}^2,\dots,-\zeta_{q_k}^{q_k-1} \quad \mbox{if } t \mbox{ is odd,}$$
or 
$$\zeta_{q_1},\zeta_{q_1}^2,\dots,\zeta_{q_1}^{q_1-1},\dots,\zeta_{q_k},\zeta_{q_k}^2,\dots,\zeta_{q_k}^{q_k-1} \quad \mbox{if } t \mbox{ is even}.$$
\end{thm}

The main ingredient of the proof of Theorem~\ref{eigen:tame} is the theory of clusters developed in \cite{DDMM}.

\begin{dfn}
Let $f(x)\in \O_F[x]$ be a squarefree monic polynomial and let $R$ be its set of roots in $\overline{F}$. 
A \emph{cluster} $\s \subseteq R$ is a non-empty set of roots of $f(x)$ of the form $\s=R\cap \mathcal{D}$ for a disc $\mathcal{D}\subseteq\overline{F}$ with respect to the $p$-adic topology. 

For a cluster $\s$ with $|\s|\geq 2$ define:
\begin{itemize}
 \item $d_\s=\min\{v(r-r'): r,r'\in \s\}$;
 \item $\s_{0}$ to be the set of maximal subclusters of $\s$ of odd size;
 \item $I_\s=\Stab_{I_F}(\s)$;
 \item $\mu_\s=\sum_{r\in R\setminus\s}v(r-r_0)$ for any $r_0\in \s$;
 \item $\lambda_\s=\frac{1}{2}(\mu_\s+d_\s |\s_0|)$;
 \item $\epsilon_\s=\begin{cases}
 \text{trivial character } \triv \text{ of }I_\s & \text{if } |\s| \text{ is even and } \ord_2 (\mu_\s \cdot  |I_F/I_\s|) \geq 1\\
 \text{order two character of }I_\s & \text{if } |\s| \text{ is even and } \ord_2 (\mu_\s \cdot  |I_F/I_\s|) < 1\\
 \text{zero representation of }I_\s & \text{otherwise };\end{cases}$
 \item $\gamma_\s:I_\s \to \C^\ast$ any character of order equal to the prime-to-$p$-part of the denominator of $|I_F/I_\s|\cdot \lambda_\s$ (with $\gamma_\s=\triv$ if $\lambda_\s=0$);
 \item $V_\s=\gamma_\s\otimes(\C[\s_0]\circleddash\triv)\circleddash \epsilon_\s$, an $I_\s$-representation.
\end{itemize}
\end{dfn}

\begin{thm}[{\cite[Theorem~$1.19$]{DDMM}}]
\label{cluster:theorem}
Let $\ell$ be a prime different from $p$, then 
$$H^1_{\acute{e}t}(C/\overline{F}, \Q_\ell)\cong H^1_{ab}\oplus (H^1_t\otimes\Sp(2))$$
as  $I_\s$-representations, with 
$$H^1_{ab}=\bigoplus_{\s\in X/I_{F}} \Ind_{I_\s}^{I_{F}}V_\s,\quad\quad H^1_t=\bigoplus_{\s\in X/I_{F}} (\Ind_{I_\s}^{I_{F}}\epsilon_\s) \circleddash \epsilon_R,$$
where $X$ is the set of clusters that are neither singletons nor (proper) disjoint unions of even size clusters, and where $\Sp(2)$ denotes the $2$\--dimensional special $\ell$-adic representation.\footnote{
Let $\ell$ and $p$ be distinct prime numbers. The \emph{special representation} $\Sp(2)$ over a local field $F/\Q_p$ is the (tame) $2$\--dimensional $\ell$-adic representation given by:
$$\Sp(2)(\tau)= \left(\begin{matrix} 1 &  t_\ell(\tau) \\ 0& 1 \end{matrix}\right), \quad 
\Sp(2)(\Frob_{\overline{F}/F})= \left(\begin{matrix} 1 & 0  \\0 & {\frac{1}{q}}\end{matrix}\right),$$
where $\tau\in I_F$, the character $t_\ell(\tau)$ is an $\ell$\--adic tame character, $\Frob_{\overline{F}/F}$ is a fixed Frobenius element and $q=\#\F$.}
\end{thm}

\begin{rmk}
When $J$ is semistable we will refer to the dimension of $H^1_t$ as the toric dimension of $J$. If the dimension of $H^1_t$ is equal to the genus $g$, we say that the reduction is totally toric.
\end{rmk}

\begin{lem}
\label{roots}
Let $f(x)\in \O_F[x]$ be a $t$-Eisenstein polynomial of degree $n$, with $(n,tp)=1$.
Then $I_F$ acts tamely on the roots of $f(x)$ and permutes them cyclically and transitively. Moreover, $v(r-r')=\frac{t}{n}$ for any two roots $r\neq r'$ of $f(x)$.
\end{lem}

\proof
Let $r$ be a root of $f(x)$. The Newton polygon of $f(x)$ has a unique slope equal to $-\frac{t}{n}$, so all the roots of $f(x)$ have valuation $\frac{t}{n}$. 
In particular, $f(x)$ is irreducible and the field $F^{\mathrm{nr}}(r)$ is a tamely ramified extension of degree $n$ of $F^{\mathrm{nr}}$. 
By uniqueness, it is Galois and its Galois group is $C_n$. As $f(x)$ is irreducible over $F^{\mathrm{nr}}$, the cyclic group $C_n$ acts transitively on the roots of $f(x)$.

Since the extension $F^{\mathrm{nr}}(r)/F^{\mathrm{nr}}$ is tame, the standard homomorphism
$$\Gal(F^{\mathrm{nr}}(r)/F^{\mathrm{nr}})\rightarrow\overline{\F}^\ast;\quad\quad \sigma \mapsto \frac{\sigma(r)}{r}$$
is injective. In particular if $\sigma$ is non-trivial then $\frac{\sigma(r)}{r}\neq 0,1$ in $\overline{\F}$, and hence $v(\sigma(r)-r)=v(r)=\frac{t}{n}$. 
As $I_F$ acts transitively on the roots, this shows that  $v(r-r')=\frac{t}{n}$ for any distinct pair of roots $r,r'$ of $f(x)$.
\endproof

\begin{lem}
\label{cluster:quantities}
Suppose that $f(x)\in \O_F[x]$ has type $t-\{q_1,{\dots},q_k\}$ with $(q_i, tp)=1$ for all $i$ and let 
$f(x)=h(x)\prod_{i=1}^k g_i(x-\alpha_i)$ be the corresponding factorisation as in Definition~\ref{type}.
Let $\beta_{i,1},\dots,\beta_{i,q_i}$ be the roots of $g_i(x-\alpha_i)$ and let $\beta_{0,1},\dots,\beta_{0,\deg h}$ be the roots of $h(x)$.
\begin{enumerate}[(i)]
 \item If $i\neq j$, then $v(\beta_{i,a}-\beta_{j,b})=0$ for all $a,b$.
 \item If $i\neq 0$, then $v(\beta_{i,a}-\beta_{i,b})=\frac{t}{q_i}$ for all $a\neq b$.
 \item $v(\beta_{0,a}-\beta_{0,b})=0$ for all $a,b$.
 \item The clusters of $f(x)$ are the whole set of roots $R$, sets $\{\beta_{i,1},\dots,\beta_{i,q_i}\}$ for every $i\neq 0$ and singleton roots.
 \item For $\s=R$, the inertia subgroup $I_F$ acts trivially on $\C[\s_0]$ and $\gamma_\s=\epsilon_\s=\triv$.
 \item For $\s=\{\beta_{i,1},\dots,\beta_{i,q_i}\}$ with $i\neq 0$, the inertia subgroup $I_F$ acts tamely on $\C[\s_0]$ and the eigenvalues of a generator $\tau$ of tame inertia are precisely the $q_i$-th roots of unity.
 Moreover, $\gamma_\s$ and $\epsilon_\s$ are also tame and 
 $$\gamma_\s(\tau)=\begin{cases}1 \mbox{ if } t \mbox{ is even}\\-1 \mbox{ if } t \mbox{ is odd,}\end{cases}
\quad \quad \epsilon_\s(\tau)=\begin{cases}0 \mbox{ if } q_i \neq 2\\ \triv \mbox{ if } q_i=2.\end{cases}$$
\end{enumerate}
\end{lem}

\proof 
\noindent
(i) The roots of $g_i(x)$ all reduce to $0$ in $\overline{\F}$, so those of $g_i(x-\alpha_i)$ all reduce to $\overline{\alpha_i}$. 
The result follows as  $\overline{\alpha_i}\not\equiv\overline{\alpha_j}$ for $i\neq j$ and  $\overline{\alpha_i}\not\equiv\overline{\beta}_{0,j}$ for all $i, j$. 

\noindent
(ii) This follows from Lemma~\ref{roots}.

\noindent
(iii) The statement follows from the definition of type $t-\{q_1,{\dots},q_k\}$ and of $h(x)$.

\noindent
(iv) Clear from (i), (ii) and (iii).

\noindent
(v) $I_F$ acts trivially on the roots of $h(x)$ and on $\{\alpha_1,\dots,\alpha_k\}$, and hence trivially on $\C[\s_0]$. 
Here $\mu_\s=d_\s=\lambda_\s=0$, so $\gamma_\s=\epsilon_\s=\triv$.

\noindent
(vi) The result for $\C[\s_0]=\C[\s]$ follows from Lemma~\ref{roots}. By (ii) $d_\s=\frac{t}{q_i}$, $\mu_\s=0$ so 
$\lambda_\s=\frac{1}{2}(\mu_\s+d_\s|\s_0|)=\frac{t}{2}$. Therefore $\gamma_\s$ is either trivial or it has order $2$ depending on the parity of $t$.
Since $\mu_\s=0$, $\epsilon_\s$ is tame and $$\epsilon_\s(\tau)=\begin{cases}0 \mbox{ if } q_i \neq 2\\ \triv \mbox{ if } q_i=2.\end{cases}$$
\endproof

\proof[Proof of Theorem~\ref{eigen:tame}]
By Lemma~\ref{cluster:quantities}~(iv) the set $X$ of clusters of $f(x)$ which are not singletons nor unions of even clusters consists of the whole set of roots $R$ and $\{\beta_{i,1},\dots,\beta_{i,q_i}\}$ for every $i\neq 0$, with $\beta_{i,j}$ as in Lemma~\ref{cluster:quantities}. 

The inertia group $I_F$ does not permute the clusters so by Theorem~\ref{cluster:theorem} we have $H^1_{\acute{e}t}(C/\overline{F}, \Q_\ell)\cong H^1_{ab}\oplus (H^1_t\otimes\Sp(2))$ with 
$H^1_{ab}=\bigoplus_{\s\in X} V_\s$ and $H^1_t=(\bigoplus_{\s\in X} \epsilon_\s) \circleddash \epsilon_R=0$, where the last equality follows from Lemma~\ref{cluster:quantities}~(vi) since $q_i\neq2$ for all $i$.

By Lemma~\ref{cluster:quantities}~(v) and~(vi), $V_\s$ is tame for each cluster $\s$. 

For $\s=R$, by Lemma~\ref{cluster:quantities}~(v) inertia acts trivially on $V_R$.

For $\s=\{\beta_{i,1},\dots,\beta_{i,q_i}\}$,  by Lemma~\ref{cluster:quantities}~(vi) the eigenvalues of $\tau$ on $V_\s$ are 
$\zeta_{q_i},\zeta_{q_i}^2,\dots,\zeta_{q_i}^{q_i-1}$ or $-\zeta_{q_i},-\zeta_{q_i}^2,\dots,-\zeta_{q_i}^{q_i-1}$ depending on whether $t$ is even or odd respectively.
In particular, $\tau$ acts semisimply on $H^1_{\acute{e}t}(C/\overline{F}, \Q_\ell)$ by an element of order dividing $2\prod_i q_i$. This proves the claim about the action of inertia on $H^1_{\acute{e}t}(C/\overline{F}, \Q_\ell)$.

As $\ell\neq p$, $H^1_{\acute{e}t}(C/\overline{F}, \Q_\ell)$ is the dual of  $T_\ell(J)\otimes \Q_\ell$ as an $I_F$-representation.
In particular, the action of $I_F$ on $T_\ell(J)$ factors through the same tame quotient and $\tau$ has the same set of eigenvalues. The result for $J[\ell]$ follows by reducing the characteristic polynomial of $\tau$ modulo $\ell$.
\endproof

\subsection{$J[\ell]$ when $\ell=p$: fundamental characters}
$\,$\\

Given an abelian variety $A/F$ with semistable reduction, a result due to Raynaud allows us to recover the eigenvalues of a generator of the tame inertia group acting on $A[p]$.
Recall that for an integer $d$ coprime to $p$, we write $\zeta_{d}$ for a primitive $d$-th root of unity, chosen such that for all divisors $d'$ of $d$ we have $\zeta_{d}^{d'}=\zeta_{\frac{d}{d'}}$.

\begin{thm}
\label{eigen:l=p}
Let $A/F$ be an abelian variety with semistable reduction. 
Then the eigenvalues of a generator of the tame inertia group on $A[p]$ are all of the form 
$$\zeta_{\substack{\; \\ p^n{-}1}}^{\substack{\sum_{i=0}^{n{-}1}a_i p^i\\ \;}}$$
\noindent
for $1\leq n\leq 2 \dim(A)$ and $0\leq a_i\leq e_F$, and where the $\zeta_{d}$ form some compatible system of roots of unity. 
\end{thm}          

For ease of reading, we will recall briefly the theory of fundamental characters. For further details see \cite[\S 1]{serre72}.

Let $I_t$ denote the tame inertia quotient of $I_F$.

A surjective homomorphism $\psi_n: I_t \to \F_{p^n}^\times$ defined by 
$$\psi_n(\sigma)=\frac{\sigma(\pi_n)}{\pi_n} \bmod{\pi_n}\quad \mbox{ where } \pi_n=\sqrt[p^n-1]{\pi_F},$$ is a fundamental character of level $n$.
The set of fundamental characters of level $n$ is the set of the $n$ characters $\psi_n$, $\psi_n^p, \ldots, \psi_n^{p^{n-1}}$; this set is independent of the choice of $\sqrt[p^n-1]{\pi_F}$.
The fundamental characters of level $n$ satisfy compatibility relations with fundamental characters of level $m$ for any integer $m$ dividing $n$:
$$\psi_n(\tau)^{ \substack{\frac{p^n-1}{p^m-1} \\ \; \\}}=\psi_m(\tau),\quad\mbox{for all}\; \tau\in I_t.$$

\proof[Proof of Theorem~\ref{eigen:l=p}]
Let $V$ be a Jordan-H\"{o}lder factor of dimension $n$ over $\F_p$ of the $I_{F}$-module $A[p]$. 
Then, by \cite[Corollaire~3.4.4]{Raynaud}, $V$ has the structure of a $1$-dimensional $\F_{p^n}$-vector space with the action of $I_{F}$ given by a character 
$\varpi \; \colon \; I_F \rightarrow \F_{p^n}^\times$, where 
$$\varpi=\psi_n^{\sum_{i=0}^{n-1} a_i p^i}$$ 
with $0\leq a_i\leq e_F$.

Let $\tau$ be a fixed generator of tame inertia. Then $\psi_n(\tau)=\zeta_{p^n{-}1}\in\F_{p^n}^\times$ and $\varpi(\tau)$ acts as multiplication by 
$\zeta_{p^n{-}1}^{\sum_{i=0}^{n{-}1}a_i p^i}$ on $\F_{p^n}$.

Let $\Phi$ be the minimal polynomial of $\zeta_{p^n{-}1}$ over $\F_p$. Since $\F_p[\zeta_{p^n{-}1}]\cong \F_{p^n}$, the minimal polynomial $\Phi$ has degree $n$ and hence its roots are 
$$\zeta_{p^n{-}1},\zeta_{p^n{-}1}^p,\zeta_{p^n{-}1}^{p^2},\dots,\zeta_{p^n{-}1}^{p^{n{-}1}}.$$
Therefore, by the Cayley\--Hamilton theorem, these are precisely the eigenvalues of multiplication by $\zeta_{p^n{-}1}$ on $V$. 

Hence, the eigenvalues of $\varpi(\tau)$ are $\zeta_{p^n{-}1}^{\sum_{i=0}^{n{-}1}a_i p^i}$.
\endproof

\subsection{Creating a transvection}
$\,$\\

Finally, we will need a criterion to ensure that some element of  $\Gal(\Q(J[\ell])/\Q)$ acts as a transvection on  $J[\ell]$. 
We again use inertia groups for achieving this.
\begin{dfn} 
Recall that a \emph{transvection} in $\GSp_{2g}(\F_\ell)$ is a unipotent element $\sigma$ such $\sigma-\Id_{2g\times 2g}$ has rank $1$.
\end{dfn}

\begin{lem}\label{inertia:transvection}
Suppose that $p\nmid 2\ell$ and $f(x)\in \O_F[x]$ has type $1-\{2\}$. Then some element of $I_F$ acts as a transvection on $J\left[\ell\right]$.
\end{lem}
\proof The model of the curve consisting of the chart $y^2=f(x)$ and the usual chart at infinity is a regular proper semistable model of $C$. 
The dual graph of the special fibre is a vertex with a loop. The homology group of the dual graph is $\Z$ with intersection pairing $(1)$, so the Tamagawa number of the Jacobian over $K^{nr}$ is $c(J/K^{nr})=\det(1)=1$ (see e.g. \cite[Theorem~3.5 and Theorem~3.8]{Papikian_2013}).

On the other hand, for a principally polarised $g$-dimensional semistable abelian variety $A$ of toric dimension $d$, the inertia group acts on $T_\ell(A)$ by block matrices of the form
$$\sigma \mapsto          
\left( \begin{array}{ccc}
\Id_d & 0 & t_\ell(\sigma)N\\
0 & \Id_{2g-2d} & 0\\
0 & 0 & \Id_d \end{array} \right)
$$
where $t_\ell$ is the $\ell$-adic tame character and $N$ is a $d\times d$ symmetric integer-valued matrix that satisfies $c(A/K^{nr})=|\mathrm{coker}(N)|$
(it is the matrix induced by the monodromy pairing composed with the principal polarisation on $A$); see e.g. \cite[\S~9,\S~10]{Grothendieck1972} or the summary in \cite[\S~3.5.1]{Dokchitser2009}. In our case $d=1$ and $c(J/K^{nr})=1$,
so $N$ is a $1\times 1$ matrix with entry $1$. In particular, picking $\sigma$ appropriately gives an element of the inertia group that acts on $J[\ell]$ as a transvection. \endproof

\section{Irreducibility}\label{section_irred}
The aim of this section is to provide explicit criteria on $f(x)$ that force irreducibility of the mod~$\ell$ Galois representation. 
The key idea is to ensure that images of local Galois groups are sufficiently large and can be patched together to guarantee global irreducibility.

\subsection{Local representations}

\begin{prop}
\label{decomposition} 
Let $C:y^2=f(x)$ be a hyperelliptic curve over a local field $F$ of odd residue characteristic $p$,  with $f(x)\in \O_F[x]$ monic and squarefree, and let $J=\Jac(C)$. 
Suppose that $f(x)$ has type $t-\{q_1,{\dots},q_k\}$ where $q_1, \dots, q_k$ are odd primes, coprime to $t$. 
Suppose moreover that the size of the residue field $\#\F$ is a primitive root modulo each of the $q_i$. Then for every prime $\ell \neq p,q_1,{\dots}, q_k$, the semisimple representation $(J[\ell]\otimes_{\F_\ell}\Flb)_{ss}$ 
decomposes as a direct sum of one \linebreak${(q_1{-}1)}$-di\-men\-sional, one ${(q_2{-}1)}$-dimensional, ${\dots}$, one ${(q_k{-}1)}$-dimensional irreduci\-ble $G_F$-subrepresentation, and all other irreducible constituents being $1$-di\-men\-sional.
\end{prop}

\proof By Theorem~\ref{eigen:tame}, $G_F$ acts tamely on  $J[\ell]$ and the non-trivial eigenvalues of a generator of tame inertia are $\pm\zeta_{q_i}^{r_{i,j}}$ 
for $r_{i,j}=1, {\dots} ,q_i{-}1$, where each sign is $+$ if $t$ is even and $-$ if $t$ is odd. 
We claim that the conclusion of the proposition holds for any semisimple $\overline{\F}_\ell$-representation $V$ with this property.

The action on $V$ factors through a finite group $G=\langle \tau,\phi\rangle$, where $\langle \tau\rangle\vartriangleleft G$ is the (tame) inertia subgroup and $\phi$ is any lift of Frobenius. 

Write $V_z$ for the $z$-eigenspace of $\tau$ on $V$.
Since $\phi \tau \phi^{-1}=\tau^{\#\F}$, and hence $\tau\phi^{-1}=\phi^{-1}\tau^{\#\F}$, it follows that $\phi^{-1}$ maps $V_z$ to $V_{z^{\#\F}}$. In particular, $\phi^{-(q_i-1)}$  is an endomorphism of $V_{\pm\zeta_{q_i}}$.

Pick $v\in V_{\pm\zeta_{q_i}}$ which is an eigenvector for the action of $\phi^{-(q_i-1)}$ on $V_{\pm\zeta_{q_i}}$ and consider the subspace 
$W=\langle v, \phi^{-1}v, {\dots}, \phi^{-(q_i-2)}v\rangle$.
Since $W$ is closed under $\tau$ and $\phi^{-1}$, it is a $G_F$-submodule of $V$.

Moreover, as $\#\F$ is a primitive root modulo $q_i$, it follows that the eigenvalues of $\tau$ on $\langle v\rangle,\dots,\langle\phi^{-(q_i-2)}v\rangle$
are precisely the non-trivial $q_i$-th roots of unity, or their negatives. In particular, as these $\tau$-eigenvalues are distinct, 
any $G$-submodule of $W$ must be a direct sum of some of the $\langle\phi^{j}v\rangle$'s. As $\phi^{-1}$ permutes these $\tau$-eigenspaces transitively, it therefore follows that $W$ is irreducible.

Now the result follows by substituting $V$ by $V/W$ and then proceeding by induction on the dimension.
\endproof

\subsection{Global representations}

\begin{lem}
\label{5out}
Let $C:y^2=f(x)$ be a hyperelliptic curve over a number field  $K$, where $f(x)\in \O_K[x]$ is a monic squarefree polynomial of degree $2g+2$. 
Suppose that $f(x)$ has type $t-\{q_1,q_2\}$ at $\p_2$ and type $t'-\{q_3\}$ at $\p_3$, 
where:
\begin{itemize}
 \item $q_1$, $q_2$ and $q_3$ are primes with $q_1\leq q_2<q_3<2g+2$ and $q_1+q_2=2g+2$;
 \item $\p_2,\p_3\nmid 2$;
 \item $t$ is coprime to $q_1q_2$, and $t'$ is coprime to $q_3$;
 \item $\#\F_{\p_2}$  is a primitive root modulo $q_1$ and $q_2$;
 \item $\#\F_{\p_3}$  is a primitive root modulo $q_3$.
\end{itemize}
Then for every prime $\ell \nmid q_1,q_2,q_3,\#\F_{\p_2}, \#\F_{\p_3}$, the $G_K$-module $J[\ell]$ is absolutely irreducible, where $J$ is the Jacobian of $C$.
\end{lem}
\proof By Proposition~\ref{decomposition}, the restriction of $J[\ell]\otimes_{\F_\ell}\Flb$ to $G_{K_{\p_3}}$ contains an irreducible $(q_3-1)$-dimensional subquotient.
Also its restriction to $G_{K_{\p_2}}$ has exactly two Jordan-Holder factors and these have dimension $(q_1-1)$ and $(q_2-1)$. 
It follows that, on the one hand, $J[\ell]\otimes_{\F_\ell}\Flb$ can have at most two Jordan-Holder factors, in which case they have dimensions $(q_1-1)$ and $(q_2-1)$, and, on the other hand, $J[\ell]\otimes_{\F_\ell}\Flb$ has a Jordan-Holder factor of dimension at least $(q_3-1)$. Hence $J[\ell]\otimes_{\F_\ell}\Flb$ is irreducible.
\endproof

\begin{thm}
\label{irred} Let $C:y^2=f(x)$ be a hyperelliptic curve over a number field  $K$, where $f(x)\in \O_K[x]$ is a monic squarefree polynomial of degree $2g+2$. 
Suppose that $f(x)$ has type $t-\{q_1,q_2\}$ at $\p_2$, type $t'-\{q_4,q_5\}$ at $\p'_2$, type $t''-\{q_3\}$ at $\p_3$, and type $t'''-\{q_5\}$ at $\p'_3$, where:
\begin{itemize}
 \item $q_1,q_2,q_3,q_4$ and $q_5$ are primes such that:
$$2g+2=q_1+q_2=q_4+q_5 \quad\quad q_4<q_1\leq q_2<q_5<q_3< 2g+2.$$
 \item $\p_2,\p'_2,\p_3,\p_3'\nmid 2 \prod_{i}q_i$ and they have distinct residue characteristics;
 \item $t$ is coprime to $q_1q_2$; $t'$ is coprime to $q_4q_5$; $t''$ is coprime to $q_3$ and $t'''$ is coprime to $q_5$;
 \item $\#\F_{\p_2}$  is a primitive root modulo $q_1$ and $q_2$;
 \item $\#\F_{\p'_2}$ is a primitive root modulo $q_4$ and $q_5$;
 \item $\#\F_{\p_3}$  is a primitive root modulo $q_3$;
 \item $\#\F_{\p'_3}$ is a primitive root modulo $q_5$.
\end{itemize}
Then for every prime $\ell$, the $G_K$-module $J[\ell]$ is absolutely irreducible, where $J$ is the Jacobian of $C$.
\end{thm}

\proof
Applying Lemma~\ref{5out} with $q_1,q_2,q_3, \p_2,\p_3$ proves the claim for all $\ell$ with $\p_2, \p_3, q_1, q_2, q_3\nmid \ell$. 
Applying the lemma again with $q_4,q_5,q_3,  \p'_2,\p_3$ proves it for all $\ell$ with $\p_3,q_3\nmid \ell$. 
By assumption $q_5>q_1,q_2$, so applying the lemma with $q_1,q_2,q_5, \p_2,\p'_3$ proves the result for $\p_3,q_3\mid \ell$.
\endproof

\section{Primitivity}\label{section_prim}
In this section $K$ is a number field and $C\colon y^2=f(x)$  a hyperelliptic curve over $K$, where 
$f(x)\in\O_K[x]$ is a monic squarefree polynomial of degree $2g+2$. As before, $J=\Jac(C)$. In this section moreover $\p$ will denote a prime of $K$.

\begin{dfn}
Let $V$ be a symplectic vector space over a field, and let $G$ be a subgroup of $\GSp (V)$. 
We say that $\{V_1,\dots,V_k\}$ is a \emph{non-trivial $G$-stable decomposition of $V$ into symplectic subspaces} if  the $V_i$ are 
 proper symplectic subspaces $V_i \subset V$,  the symplectic pairing is non-degenerate on $V_i$, and there is a homomorphism $\phi : G \rightarrow S_k$ such that $V=\oplus_{i=1}^k V_i$ and $\sigma(V_i)=V_{\phi(\sigma)(i)}$ for $\sigma \in G$.
\end{dfn}

\begin{dfn}
\label{imprimitive_def}
Let $V$ be a symplectic vector space over a field, and let $G$ be a subgroup of $\GSp (V)$.
Suppose that $V$ has no proper $G$-stable subspace. Recall that $V$ is an \emph{imprimitive} $G$-module 
if there is no non-trivial $G$-stable decomposition of $V$ into symplectic subspaces. If $V$ is not an imprimitive $G$-module, then it is a \emph{primitive} $G$-module.
\end{dfn}

\subsection{Quasi-unramified representations}

\begin{dfn}[Quasi-unramified representation]
We will say that a symplectic $\F_\ell$-representation $V$ of $G_K$ is \emph{quasi-unramified} if for every $G_K$-stable decomposition 
$V=\oplus_{i=1}^k V_i$ into symplectic $\F_\ell$-subspaces, the permutation action of $G_K$ on $\{V_1,\dots,V_k\}$ is unramified at every prime of $K$. 
We will say that $V$ is \emph{strongly quasi-unramified} if the same condition holds for decompositions of
$V{\otimes_{\F_\ell}}\overline{\F}_\ell$ into symplectic $\overline{\F}_\ell$-subspaces.
\end{dfn}

Note that strongly quasi-unramified implies quasi-unramified.

\begin{prop}
\label{quasiun:primitive}
Let $K$ be a number field which does not have everywhere unramified extensions. 
If $V$ is an irreducible quasi-unramified symplectic representation of $G_K$, then $V$ is primitive.
\end{prop}

\proof Suppose $V$ admits a $G_K$-stable decomposition into symplectic subspaces. 
Since $V$ is irreducible, the associated homomorphism $\phi$ is transitive and, in particular, non-trivial. 
By definition of quasi-unramified symplectic representation, the kernel of $\phi$ cuts out a proper unramified extension of $K$. 
Hence, 
$V$ is primitive. 
\endproof

\begin{rmk}
By the Hermite\--Minkowski theorem, $\Q$ satisfies the hypotheses of Proposition~\ref{quasiun:primitive}. Similarly, the same holds for $\Q(\zeta_3)$.
\end{rmk}

\subsection{Criteria for being quasi-unramified }

\begin{dfn}[Admissible polynomials]\label{adm_poly}
We will say that $f(x)\in \O_K[x]$ is \linebreak{\emph{$\ell$-admissible at $\p$}} if for every $G_K$-stable decomposition 
$J[\ell]\otimes\overline{\F}_\ell=\oplus_{i=1}^k V_i$ into symplectic $\overline{\F}_\ell$-subspaces, $I_{\p}$ acts trivially on $\{V_1,\dots,V_k\}$.

We will say that $f(x)\in \O_K[x]$ is \emph{admissible at $\p$} if it is $\ell$-admissible at $\p$ for every odd prime number $\ell$ not divisible by $\p$.
\end{dfn}
\begin{prop}
\label{strong}
Let $\ell$ be an odd prime number and suppose that:
\begin{enumerate}[(1)]
 \item $f(x)$ is admissible at all $\p\nmid \ell$;
 \item $f(x)$ is $\ell$-admissible at all $\p\mid\ell$.
\end{enumerate}
Then $J[\ell]$ is strongly quasi-unramified.
\end{prop}

\proof Direct from the definition. \endproof

The following criterion is another way of ensuring that $J[\ell]$ is quasi-unramified for certain primes $\ell$:
\begin{prop}
\label{nonsemistable}
If $J[\ell]$ is irreducible and there exists a prime number $q$ such that 
\begin{itemize}
 \item $g+1<q<2g+2$,
 \item $\ell$ is a primitive root modulo $q$,
 \item $f$ has type $t-\{q\}$ at $\p$ for some $\p\nmid 2\ell q t$,
\end{itemize}
then $J[\ell]$ is quasi-unramified.
\end{prop}

\proof Suppose that $J[\ell]=\oplus_{i=1}^k V_i$ is a non-trivial $G_K$-stable decomposition into symplectic $\F_\ell$-subspaces. Since $J[\ell]$ is irreducible and the $V_i$ are symplectic subspaces, $\dim V_i = \dim V_j\geq 2$. In particular, $k< g+1 <q$. 
By Theorem~\ref{eigen:tame}, $I_\p$ acts on $J[\ell]$ through a cyclic quotient $\langle \tau \rangle \cong C_{q}$ or $C_{2q}$, with 
$\tau^2$ having eigenvalues all the primitive  $q$-th roots of unity, and all other eigenvalues $1$. 
By hypothesis $k<q$, so $\tau^2$ has to preserve each of the $V_i$. Without loss of generality $\tau^2:V_1 \to V_1$ has $\zeta_{q}$ as an eigenvalue.
The minimal polynomial of $\zeta_{q}$ over $\F_\ell$ has degree $q-1$, so $\frac{2g}{k}=\dim V_1\geq q-1>g$, and so $k=1$.
\endproof

In the rest of this section we will give criteria for $f(x)$ to satisfy the hypotheses of Proposition~\ref{strong}.

\subsection{Admissible polynomials}

\begin{lem}
\label{inertia:noperm}
If $J$ is semistable at $\p$ then $f(x)$ is admissible at $\p$. 
\end{lem}

\proof
Let $\ell$ be an odd prime with $\p\nmid \ell$. 
Suppose that $J[\ell]\otimes\overline{\F}_\ell=\oplus_{i=1}^k V_i$ is a non-trivial $G_K$-stable decomposition of $J[\ell]\otimes\overline{\F}_\ell$ into symplectic $\overline{\F}_\ell$-subspaces.

By \cite[Corollaire~3.5.2]{Grothendieck1972}, $I_{\p}$ acts unipotently on $J[\ell]$ with $(\sigma-1)^2=0$ for all $\sigma\in I_{\p}$. 
By Lemma~\ref{mat:wild}, every $\sigma\in I_{\p}$ fixes each of the $V_i$ and so $f(x)$ is $\ell$-admissible at $\p$. 
\endproof

\begin{lem}
\label{inertia:adm_q1q2}
Let $\p$ be a prime of odd residue characteristic $p$. 
If $f(x)\in \O_K[x]$ has type $t-\{q_1,q_2\}$ at $\p$ for $t$ odd and odd primes $q_1,q_2$ different from $p$, with $q_1+q_2=2g+2$, then $f(x)$ is admissible at $\p$.
\end{lem}

\proof Let $\ell\neq 2,p$ be a prime. Suppose that $J[\ell]\otimes\overline{\F}_\ell=\oplus_{i=1}^k V_i$ is a non-trivial $G_K$-stable decomposition into symplectic $\overline{\F}_\ell$-subspaces.
By Theorem~\ref{eigen:tame}, $I_{\p}$ acts tamely on $J[\ell]$ through a cyclic quotient of order dividing $2q_1q_2$. 
Let $\tau$ be a fixed generator of tame inertia. We need to show that $\tau$ acts trivially on $\{V_1,\dots,V_k\}$. 
Again by Theorem~\ref{eigen:tame}, $\tau$ has eigenvalues $-\zeta_{q_1},\dots,-\zeta_{q_1}^{q_1-1},-\zeta_{q_2},\dots,-\zeta_{q_2}^{q_2-1}$.

If $\ell\neq q_1,q_2$, then no subset of the eigenvalues is closed under multiplication by either $-1$,  $\zeta_{q_1}$ or $\zeta_{q_2}$, and so by Lemma~\ref{lem:cyclic}, $\tau$ cannot permute $\{V_1,\dots,V_k\}$.

If $\ell=q_1\neq q_2$, then by the same argument $\tau$ cannot have an orbit on $\{V_1,\dots,V_k\}$ of length $2$.
Moreover, $\tau$ does not have an eigenvalue of multiplicity $q_1$, so by Lemma~\ref{lem:cyclic} $\tau$ cannot have an orbit of lenght divisible by $q_1$.
Furthermore, no set of $2q_2$ eigenvalues is closed under multiplication by $\zeta_{q_2}$ and the $V_i$ are symplectic (even dimension), so 
$\tau$ cannot have an orbit of length divisible by $q_2$ either. 

Finally, if $\ell=q_1=q_2=g+1$, then $\tau$ cannot cyclically permute $q_1$ symplectic subspaces since $\dim_{\F_\ell}J[\ell]=2g<2q_1$.
It also cannot have an orbit of length $2$ as no subset of the eigenvalues is closed under multiplication by ${-}1$. 
\endproof

\begin{lem}
\label{inertia:adm_q3}
Let $\p$ be a prime of odd residue characteristic $p$. 
If $f(x)\in \O_K[x]$ has type $2{-}\{q\}$ at $\p$ where $q$ is an odd prime $g+1<q<2g+2$, then $f(x)$ is admissible at $\p$.
\end{lem}

\proof
Let $\ell\neq 2,p$ be a prime. By Theorem~\ref{eigen:tame}, $I_{\p}$ acts tamely on $J[\ell]$ through a cyclic quotient of order dividing $2q$, and with a 
generator of tame inertia $\tau$ having non-trivial eigenvalues $\zeta_{q},\dots,\zeta_{q}^{q-1}$, each with multiplicity $1$ (unless $\ell=q$ in which case all eigenvalues are $+1$).

Since $\ell\neq 2$, the order of the image of $I_\p$ is $q$ or $1$. Clearly, since $q>g+1$, inertia cannot permute $q$ symplectic blocks.\endproof

\subsection{$p$-admissible polynomials}
$\;$\\

Let us address condition $(2)$ in Proposition~\ref{strong}.

\begin{prop}
\label{primitivelp}
If $C/K$ is semistable at $\p$ of residue characteristic $p$, with 
$$p> \max(g, 2 e_{K_\p} +1),$$ 
where $e_{K_\p}$ is the ramification degree of $K_\p/\Q_p$, then $f(x)$ is $p$-admissible at $\p$. 
\end{prop}

\proof Suppose that $J[p]\otimes\overline{\F}_p=\oplus_{i=1}^m V_i$ is a non-trivial $G_{K_\p}$-stable decomposition into symplectic subspaces. 

As $p>g$ the wild inertia group cannot permute the subspaces, as each orbit must have either size $1$ or size divisible by $p$.

By Theorem~\ref{eigen:l=p}, the eigenvalues of a (fixed) generator $\tau$ of the tame inertia group are of the form $\zeta_{p^k-1}^{\sum_{i=0}^{k-1}a_i p^i}$
for some $1\leq k\leq 2\dim J$ and $0\leq a_i\leq e_{K_\p}$. 
In particular, if $\zeta_x$ is a root of unity such that $\zeta_x^{\frac{x}{p^k-1}}=\zeta_{p^k-1}$ for all $k$, then each eigenvalue is of the form 
$$\zeta_x^{tx} \quad \quad \mbox{for some}\quad 0\leq t \leq \frac{e_{K_\p}}{p-1}<\frac{1}{2}.$$
This set has no subset closed under multiplication by $j$-th roots of unity for any $j\leq g$ (as $g < p$). 
Thus by Lemma~\ref{lem:cyclic}, $\tau$ cannot permute $\{V_1,\dots,V_m\}$.
\endproof

\begin{rmk}
The result and the proof of Proposition~\ref{primitivelp} also hold for abelian varieties. 
Let $A$ be a $g$-dimensional semistable abelian variety over a local field $F/\Q_p$. 
Suppose that $A[p]\otimes\overline{\F}_p=\oplus_{i=1}^m V_i$ is a $G_F$-stable decomposition into symplectic subspaces.
If $p> \max(g, 2e_F+1)$ then $I_{F}$ does not permute $\{V_1,\dots,V_m\}$.
\end{rmk}

\begin{prop}
\label{toric_quasi}
Let $\p$ be a prime of odd residue characteristic $p$ with 
$$e_{K_\p(\zeta_p)/K_\p} \neq 2.$$
If $J$ has totally toric reduction at $\p$ then $f(x)$ is $p$-admissible at $\p$. 
\end{prop}

\proof Suppose that $J[p]\otimes\overline{\F}_p=\oplus_{i=1}^k V_i$ is a $G_K$-stable decomposition of $J[p]$ into symplectic $\overline{\F}_p$-subspaces. 
The inertia group $I_{\p}$ acts on $T_p(J)$ as
 $$\sigma \mapsto \left(
\begin{array}{c|c}
\chi(\sigma)\Id_{g\times g} & \ast \\ \hline
{0} & {\Id_{g\times g}} 
\end{array}\right),$$
where $\chi$ is the cyclotomic character (as follows from the Raynaud parametrization $J(\overline{K}_{\p})\cong (\overline{K}_{\p}^\times)^g/$lattice, so that $0\to (\mu_{p^n})^g\to J(\overline{K}_{\p})[p^n] \to (\Z/p^n)^g \to 0$). 
In particular, the action of $I_{\p}$ on $J[p]$ factors through a group of the form 
$G=W\rtimes H$ where $W$ is a $p$-group with all matrices satisfying $(M-\Id)^2=0$, and $H=\langle \tau \rangle$ is a cyclic group of order $e_{K_\p(\zeta_p)/K_\p}$. 
The eigenvalues of $\tau$ are $1$ and $\chi(\tau)$, both with multiplicity $g$.

By Lemma~\ref{mat:wild}, $W$ acts trivially on $\{V_1,\dots,V_k\}$. Since $\chi(\tau)\neq -1$, the eigenvalues of $\tau$ have no subset closed under multiplication by any root of unity.
Thus, by Lemma~\ref{lem:cyclic} $\tau$ cannot permute any of the $V_i$ either. Therefore, $I_{\p}$ acts trivially on $\{V_1,\dots,V_k\}$, as required.

\endproof

\subsection{Miscellaneous linear algebra}

\begin{lem}\label{lem:cyclic}
Let $V=\oplus_{i=1}^{k} V_i$ be a finite dimensional vector space over a field $L$. Let $T \colon V \rightarrow V$ be an $L$-linear map
such that $T(V_i)=V_{i+1}$ (the indices considered modulo $k$). 
If the eigenvalues of $T^k$ on $V_1$ are $\alpha_1, \dots, \alpha_d$ (with multiplicity), then the eigenvalues of $T$ on $V$ (with multiplicity) are 
$$\zeta_k^j\sqrt[\leftroot{-2}\uproot{2}k]{\alpha_i}$$
for $i=1,\dots,d$ and $j=0,\dots,k-1$.
\end{lem}
\proof  Without loss of generality, suppose that $L$ is algebraically closed. Pick $v\in V$ such that $Tv=\beta v$. Write $v=\sum_{i=1}^k v_i$ for $v_i\in V_i$, so 
$Tv_i=\beta v_{i+1}$. 

On the subspace $W=\langle v_1,\dots, v_k \rangle$ the map $T^k$ acts as multiplication by $\beta^k$, so $T^k-\beta^k=0$ on $W$. 
The minimal polynomial of $T$ on $W$ must have at least degree $k$, so by the Cayley\--Hamilton theorem the characteristic polynomial of $T$ on $W$ is $x^k-\beta^k$. 
Hence its eigenvalues are $\zeta_k^j\beta$ for $j=0,\dots,k-1$. Now take $V'=V/W=\oplus_{i=1}^{k} V_i/\langle v_i \rangle$ and proceed by induction on the dimension.
\endproof

\begin{lem}
\label{mat:wild}
Let $\ell$ be an odd prime and $V$ an $\Flb$-vector space. Suppose $M:V{\to}V$ is a linear map and satisfies $(M-\Id)^2=0$. 
Then there is no set of linearly independent subspaces $V_1,\dots,V_k$ of $V$ (for $k>1$) which are cyclically permuted~by~$M$.
\end{lem}
\proof Since $0=(M-\Id)^\ell=M^\ell-\Id$, either $M=\Id$ or $M$ has order $\ell$. So, if it permutes a set of linearly independent subspaces $V_1,\dots,V_k$ cyclically,
then $k=\ell>2$. Now if $v\in V_1\setminus\{0\}$, then
$$0=(M-\Id)^2 v= M^2v-2Mv+v,$$
which gives a contradiction since $v\in V_1\setminus\{0\}$, $2Mv\in V_2$ and $M^2v\in V_3$.
\endproof

\section{Surjectivity}\label{section_surj}

\subsection{Generating symplectic groups}
$\;$\\

We make use of the following classification of subgroups of $\GSp_{2g}(\F_\ell)$ containing a transvection, due to Hall and Arias-de-Reyna, Dieulefait, Wiese.

\begin{thm}[{\cite[Theorem~1.1]{hall08}; \cite[Theorem~1.1]{disawi}}]
\label{classification}
Let $\ell \ge 5$ be a prime and let $V$ be a symplectic $\F_\ell$-vector space.
Let $G$ be a subgroup of $\GSp (V)$ such that:
\begin{enumerate}[$(i)$]
\item $G$ contains a transvection;
\item $V$ is an $\F_\ell$-irreducible $G$-module;
\item $V$ is a primitive $G$-module.
\end{enumerate}
Then $G$ contains $\Sp(V)$. The same is true for $\ell=3$, provided that $V\otimes\overline{\F}_3$ is an irreducible and primitive $G$-module.
\end{thm}

\begin{prop}
If $\ell\geq 5$ and $V$ is an irreducible quasi-unramified symplectic representation of $G_\Q$, then 
the image of $G_\Q$ contains $\Sp(V)$ provided that some element of $G_\Q$ acts as a transvection.
The same holds for $\ell=3$ provided that $V$ is also absolutely irreducible and strongly quasi-unramified.
\end{prop}

\proof By Lemma~\ref{quasiun:primitive} the representation is primitive. The result follows from Theorem~\ref{classification}.
\endproof

\subsection{Symplectic representations and abelian varieties}

\begin{thm}
\label{classification_quasi:number_fields}
Let $\ell\geq 5$ be a prime and let $A/K$ be a principally polarized abelian variety of dimension $g$ over a number field $K$. 
If the $G_K$-action on $A[\ell]$ is irreducible, primitive and contains a transvection, then the image of $G_K$ contains $\Sp_{2g}(\F_\ell)$. 
Moreover, the same holds for $\ell=3$ provided that $A[3]\otimes\overline{\F}_3$ is also irreducible.
\end{thm}

\proof The result follows directly from Theorem~\ref{classification}.
\endproof

\begin{lem}
\label{index_cycl}
Let $\ell$ be a prime and let $A/K$ be a principally polarized abelian variety of dimension $g$ over a number field $K$. 
Let  $G=\Gal(K(A[\ell])/K)$. Then $[G:G\cap \Sp_{2g}(\F_\ell)]=[\Q(\zeta_\ell): K\cap\Q(\zeta_\ell)]$.
\end{lem}
\proof
Let $t:\GSp_{2g}(\F_\ell)\to\F_\ell^\times$ be the group homomorphism which maps an element ${M\in\GSp_{2g}(\F_\ell)}$ to the corresponding multiplier through the symplectic pairing, that is the element $m \in \F_\ell^\times$ such that for all $u,v\in \F_\ell^{2g}$, $\langle Mu, Mv \rangle=m \langle u, v \rangle$.
The kernel of this homomorphism is $\Sp_{2g}(\F_\ell)$. 

Since the abelian variety is principally polarized, the symplectic pairing on $J[\ell]$ is the mod~$\ell$ Weil pairing: for all $\sigma \in G_K$ and for all $v,w \in J[\ell]$ we have $\langle \sigma v, \sigma w \rangle =\chi(\sigma) \langle v,w \rangle$. Therefore the homomorphism $t$ restricted to $G$ is the cyclotomic character and 
$$[G:G\cap \Sp_{2g}(\F_\ell)]=|\mbox{image of }\chi \mbox{ on } G|=[\Q(\zeta_\ell): K\cap\Q(\zeta_\ell)].$$\endproof

\begin{cor}
\label{cycl}
Let $\ell\geq 5$ be a prime and let $J/K$ be the Jacobian of a curve of genus $g$ over a number field $K$. 
If the $G_K$-action on $J[\ell]$ is irreducible, primitive and contains a transvection, then $\Gal(K(J[\ell])/K)$ contains $\Sp_{2g}(\F_\ell)$ with index 
$[\Q(\zeta_\ell): K\cap\Q(\zeta_\ell)]$.
\noindent
Moreover, the same holds for $\ell=3$ provided that $J[3]\otimes\overline{\F}_3$ is also irreducible.
\end{cor}
\proof Since $J$ is principally polarized (see \cite[Summary~6.11]{milne_jacobians}), Theorem~\ref{classification_quasi:number_fields} implies that the image of $G_K$ contains  $\Sp_{2g}(\F_\ell)$. 
The result follows from Lemma~\ref{index_cycl}.\endproof

\begin{thm}
\label{full_Q}
Let $\ell\geq 5$ be a prime and let $J/\Q$ be the Jacobian of a curve of genus $g$. 
If the $G_\Q$-action on $J[\ell]$ is irreducible, quasi-unramified and contains a transvection, then $\Gal(\Q(J[\ell])/\Q)\cong \GSp_{2g}(\F_\ell)$. 
Moreover, the same holds for $\ell=3$ provided that $J[3]$ is also absolutely irreducible and strongly quasi-unramified.
\end{thm}

\proof By Lemma~\ref{quasiun:primitive} the representation is primitive. The result follows from Corollary~\ref{cycl}.\endproof

\section{Maximal Galois images over $\Q$}\label{section_maximal}
We now put together the results from \S~2-\S~5 to produce hyperelliptic curves over $\Q$ with maximal Galois images.
In this section $C\colon y^2=f(x)$ denotes a hyperelliptic curve over $\Q$, where $f(x)\in\Z[x]$ is a monic squarefree polynomial of degree $2g+2$. As before, $J=\Jac(C)$. 

For the rest of the section we will refer to the following hypotheses on the genus and on $f(x)$:

\smallskip
\begin{itemize}
 \item[(G$+\epsilon$)] There exist primes $q_1,q_2$ and $q_3$ such that:
$$2g+2=q_1+q_2, \quad\quad q_1\leq q_2<q_3< 2g+2.$$
 \item[(2G$+\epsilon$)] There exist primes $q_1,q_2,q_3,q_4,q_5$ such that:
  $$2g+2=q_1+q_2=q_4+q_5, \quad\quad q_4<q_1\leq q_2<q_5<q_3< 2g+2.$$
 \item[(2T)] $f(x)$ has type $1-\{2\}$ at distinct primes $p_t,p_t'>g$.
 \item[(TT)] $J$ has totally toric reduction at all odd primes $\ell\leq g$.
 \item[($p_2$)] $f(x)$ has type $1-\{q_1,q_2\}$ at a prime $p_2>2g+2$, which is a primitive root modulo~$q_1, q_2$ and $q_3$.
 \item[($p_3$)] $f(x)$ has type $2{-}\{q_3\}$ at a prime $p_3>2g+2$, which is a primitive root modulo~$q_3$.
 \item[($p_2'$)]$f(x)$ has type $1-\{q_4,q_5\}$ at a prime $p_2'>2g+2$, which is a primitive root modulo~$q_3,q_4$ and $q_5$.
 \item[($p_3'$)] $f(x)$ has type $2{-}\{q_5\}$ at a prime $p_3'>2g+2$, which is a primitive root modulo~$q_5$.
 \item[(adm)]$f(x)$ is admissible at all primes $p$ (see Definition~\ref{adm_poly}).
 \item[(ss)] $C$ is semistable at all primes $p\notin\{p_2,p_2',p_3,p_3'\}$.
 \item[(3)]$p_2,p_3\equiv 1 \bmod{3}$.
 \item[($S_{2g+2}$)] There exist two primes $p_{irr}$ and $p_{lin}$ such that $f(x)$ modulo $p_{irr}$ is irreducible, and $f(x)$ modulo $p_{lin}$ factors as an irreducible polynomial times a linear factor.
\end{itemize}
\smallskip

Theorem~\ref{out:7} requires the Goldbach conjecture like hypothesis $(\mathrm{G}+\epsilon)$ and produces curves with maximal mod~$\ell$ Galois images at all but a small set of primes.
Theorem~\ref{alldone} requires the stronger hypothesis $(2\mathrm{G}+\epsilon)$ but guarantees maximality at all $\ell$ simultaneously. 

\begin{rmk}
Note that hypothesis $(\mathrm{adm})$ is automatically satisfied if hypotheses $(p_2)$, $(p_3)$, $(p_2'), (p_3')$ and $(\mathrm{ss})$ hold.
The polynomial $f(x)$ is admissible at all primes: Lemma~\ref{inertia:adm_q1q2} and \ref{inertia:adm_q3} ensure admissibility at $p_2,p_3,p_2'$ and $p_3'$, while 
hypothesis $(\mathrm{ss})$ and Lemma~\ref{inertia:noperm} guarantee admissibility at all other primes.
\end{rmk}

\begin{thm}
\label{out:7}
Suppose $f(x)\in \Z[x]$ satisfies $(\mathrm{G}+\epsilon)$, $(2\mathrm{T})$, $(p_2)$, $(p_3)$ and $(\mathrm{adm})$. 
Then $\Gal(\Q(J[\ell])/\Q)\cong \GSp_{2g}(\F_\ell)$ provided that $\ell \neq 2,3,q_1,q_2,q_3,p_2,p_3$ and either
\begin{enumerate}[(i)]
 \item $\ell>g$ and $J/\Q_\ell$ is semistable, or
 \item $J/\Q_\ell$ has totally toric reduction, or
 \item $\ell$ is a primitive root modulo $q_3$.
\end{enumerate}
\end{thm}

\proof 
By Lemma~\ref{inertia:transvection}, Hypothesis $(2\mathrm{T})$ ensures the existence of a transvection in $\Gal(\Q(J[\ell])/\Q)$.

By $(\mathrm{G}+\epsilon)$, $(p_2)$ and $(p_3)$, the hypotheses of Lemma~\ref{5out} are satisfied, so $J[\ell]$ is absolutely irreducible for every prime 
$\ell \neq q_1,q_2,q_3,p_2,p_3$. 

In case $(i)$, by Proposition~\ref{primitivelp} and hypothesis $(\mathrm{adm})$, the conditions of Proposition~\ref{strong} are satisfied, so $J[\ell]$ is strongly quasi-unramified.

In case $(ii)$, by Proposition~\ref{toric_quasi} and hypothesis $(\mathrm{adm})$, the conditions of Proposition~\ref{strong} are again satisfied, so $J[\ell]$ is strongly quasi-unramified.

In case $(iii)$, since $J[\ell]$ is irreducible and hypothesis $(p_3)$ holds, Proposition~\ref{nonsemistable} shows that $J[\ell]$ is quasi-unramified.

Therefore, $\Gal(\Q(J[\ell])/\Q)\cong \GSp_{2g}(\F_\ell)$ by Theorem~\ref{full_Q}.
\endproof

\begin{rmk}\label{out:5}
The theorem can be easily extended to include $\ell=p_2$ and $p_3$ by requiring that there is a second pair of primes $r_2, r_3$ satisfying the same properties as $p_2$ and $p_3$ in hypotheses $(p_2)$ and  $(p_3)$.
\end{rmk}

From Theorem~\ref{out:7} we have the following immediate corollary:
\begin{cor}
If $f(x)$ also satisfies $(\mathrm{TT})$ then $\Gal(\Q(J[\ell])/\Q)\cong \GSp_{2g}(\F_\ell)$ for every prime $\ell$, except possibly for
\begin{enumerate}[(i)]
 \item $\ell=2,3,q_1,q_2,q_3,p_2,p_3,$ and 
 \item $\ell$ where $J/\Q_\ell$ is not semistable that are not primitive generators modulo~${q_3}$.
\end{enumerate}
\end{cor}

\begin{thm}
\label{alldone}
Suppose $f(x)\in \Z[x]$ satisfies $(2\mathrm{G}+\epsilon)$, $(2\mathrm{T})$, $(p_2)$, $(p_3)$, $(p_2'), (p_3')$, $(\mathrm{TT})$ and $(\mathrm{ss})$.
Then $\Gal(\Q(J[\ell])/\Q)\cong \GSp_{2g}(\F_\ell)$ for all primes $\ell \neq 2,3$. 

Moreover, if $f(x)$ also satisfies $(3)$ then $\Gal(\Q(J[3])/\Q)\cong \GSp_{2g}(\F_3)$, and if $f(x)$ satisfies $(S_{2g+2})$ then $\Gal(\Q(J[2])/\Q)\cong S_{2g+2}$.
\end{thm}

\proof \underline{Case $\ell\geq 5$}. The hypotheses of Theorem~\ref{irred} are satisfied by $(2\mathrm{G}+\epsilon)$, $(p_2)$, $(p_3)$, $(p_2')$ and $(p_3')$, so $J[\ell]$ is absolutely irreducible for every prime $\ell$.  

By hypothesis $(2\mathrm{T})$, Lemma~\ref{inertia:transvection} ensures that for every prime $\ell$ there exists a transvection in $\Gal(\Q(J[\ell])/\Q)$.

The polynomial $f(x)$ is admissible at all primes: Lemmas~\ref{inertia:adm_q1q2} and \ref{inertia:adm_q3} ensure admissibility for $p_2,p_3,p_2'$ and $p_3'$, while 
hypothesis $(\mathrm{ss})$ and Lemma~\ref{inertia:noperm} guarantee admissibility at all other primes.

If $\ell \leq g$, then $J$ has totally toric reduction by hypothesis $(\mathrm{TT})$.
Since $\ell\geq 5$, we have that $e_{\Q_\ell(\zeta_\ell)/\Q_\ell}\neq 2$. Therefore, by Proposition~\ref{toric_quasi} $J[\ell]$ is $\ell$-admissible at $\ell$. 
If $\ell>g$ and $\ell \notin\{p_2,p_3,p_2',p_3'\}$, then $J$ is semistable at $\ell$ by hypothesis $(\mathrm{ss})$ and so by Proposition~\ref{primitivelp} $J[\ell]$ is $\ell$-admissible at $\ell$. 

If $\ell \notin\{p_2,p_3,p_2',p_3'\}$, then by Proposition~\ref{strong} $J[\ell]$ is strongly quasi-unramified. 
If $\ell \in \{p_2,p_3,p_2',p_3'\}$, then Proposition~\ref{nonsemistable} with $q=q_3$ (or $q_5$) and $p=p_3$ (or $p_3'$) shows have that $J[\ell]$ is quasi-unramified.

Therefore, by Theorem~\ref{full_Q} we have that $\Gal(\Q(J[\ell])/\Q)\cong \GSp_{2g}(\F_\ell)$.\\ 
\smallskip

\underline{Case $\ell=3$}.
We will show that $\Gal(K(J[3])/K)$ contains $\Sp_{2g}(\F_3)$ for $K=\Q(\zeta_3)$. By Lemma~\ref{index_cycl} then we have $\Gal(\Q(J[3])/\Q)\cong \GSp_{2g}(\F_3)$.

Note that hypothesis $(2\mathrm{G}+\epsilon)$ forces $g\geq6$. Since $p_t>g\geq 6$, it is unramified in $K/\Q$ and $f(x)$ has type $1-\{2\}$ at all primes dividing $p_t$. 
The existence of a transvection in $\Gal(K(J[3])/K)$ is ensured by Lemma~\ref{inertia:transvection}.

Let $\p_2,\p_3$ be primes of $K$ above $p_2$ and $p_3$ respectively. 
By hypothesis $(3)$, $p_2$ and $p_3$ split in $K$, so $\#\F_{\p_2}=p_2$ and $\#\F_{\p_3}=p_3$, and $f(x)$ has type $1-\{q_1,q_2\}$ at $\p_2$ and type $2{-}\{q_3\}$ at $\p_3$.
Let us remark that $p_2,p_3,q_1,q_2,q_3\neq 3$ since $q_4<q_1\leq q_2<q_3$ and $p_2,p_3>g\geq 6$. By Lemma~\ref{5out}, the $G_K$-module $J[3]$ is absolutely irreducible.

The hyperelliptic curve $C/K$ is semistable at all primes $\p \nmid p_2,p_2',p_3,p_3'$ by hypothesis $(\mathrm{ss})$.
As $p_2,p_2',p_3,p_3'$ are unramified in $K/\Q$, $f(x)$ has the same type above these primes in $K$ as over $\Q$. In particular $f(x)$ is admissible at all primes by 
Lemmas~\ref{inertia:noperm},~\ref{inertia:adm_q1q2} and~\ref{inertia:adm_q3}. 
By hypothesis $(\mathrm{TT})$, $J$ has totally toric reduction at $\lambda\mid 3$ and so, by Proposition~\ref{toric_quasi} $J[3]$ is $3$-admissible at $\lambda$.
By Proposition~\ref{strong} $J[\ell]$ is strongly quasi-unramified. 
By Proposition~\ref{quasiun:primitive} $J[3]\otimes\overline{\F}_3$ is primitive. 

Theorem~\ref{classification_quasi:number_fields} shows that $\Sp_{2g}(\F_3) \subseteq \Gal(K(J[3])/K)$, as required.\\ 

\smallskip

\underline{Case $\ell=2$}. 
Recall\footnote{As in \cite{Cornelissen2001}, $J[2]$ is generated by divisors $D_i=(t_i,0)-(t_1,0)$ for $i=2,\dots, 2g+2$  subject to the unique relation $\sum_i D_i=0$, where the $t_i$ are the roots of $f(x)$. Clearly $\Q(J[2])$ is contained in the splitting field of $f(x)$. Conversely, it is easy to see that if $\deg(f)>4$ and $\sigma \in \Gal(f)$ satisfies $\sigma(D_i)=D_i$ for all $i$ then $\sigma$ is trivial.} that $\Q(J[2])$ is the splitting field of $f(x)$. We just need to ensure that the Galois group of $f(x)$ is the full symmetric group $S_{2g+2}$.

Hypothesis $(S_{2g+2})$ guarantees the existence of primes $p_{irr}$ and $p_{lin}$ such that 
$f(x)$ modulo $p_{irr}$ is irreducible, and $f(x)$ modulo $p_{lin}$ factors as an irreducible polynomial times a linear factor.
These factorisations ensure the existence of a $2g+2$ cycle and a $2g+1$ cycle in the Galois group of $f(x)$. 

By hypothesis $(2\mathrm{T})$, the inertia group at $p_t$ acts as a transposition on the roots of $f(x)$. 

Since using $2g+2$ and $2g+1$ cycles it is possible to conjugate a transposition to any other transposition, and the symmetric group is generated by transpositions, 
we deduce that  the Galois group of the splitting field of $f(x)$ is $S_{2g+2}$, as required.
\endproof

\begin{rmk}
\label{exceptional_list}
Hypothesis $(2\mathrm{G}+\epsilon)$ does not hold for $g=0,1,2,3,4,5,7$ and $13$, but we expect it to hold for all other $g$, and have numerically verified it for $g\leq 10^7$.

For this exceptional list of small genera our method still makes it possible to find hyperelliptic curves with $\Gal(\Q(J[\ell])/\Q)\cong \GSp_{2g}(\F_\ell)$ for all but a small set of primes $\ell$
(see Theorem~\ref{out:7}, Remark~\ref{out:5}, hypothesis~$(S_{2g+2})$ and the proofs of cases $\ell=2,3$ of Theorem~\ref{alldone}):
\begin{center}
\begin{tabular}{ l|l }
Genus & primes excluded \\
\hline
  $2$ & $3,5$\\
  $3$ & $3,5,7$ \\
  $4$ & $5,7$ \\
  $5$ & $5,7,11$ \\
  $7$ & $5,11,13$ \\
  $13$& $11,17,23$. 
\end{tabular}
\end{center}
\end{rmk}

\section{Congruence conditions}\label{section_congruence}

\subsection{Main theorem: explicit curves}
$\;$\\

The main result of this section is the following explicit version of Theorem~\ref{alldone}:

\begin{thm}
\label{final:congruence_th}
Let $g$ be a positive integer such that there exist primes $q_1,q_2,q_3,q_4,q_5$ with 
$2g+2=q_1+q_2=q_4+q_5$ and $q_4<q_1\leq q_2<q_5<q_3<2g+2$.
Let 
$$f_0(x)=x^{2g+2}+a_{2g+1}x^{2g+1}+\dots+a_1 x+a_0 \in \Z[x]$$ 
be a polynomial such that
\begin{itemize}
 \item $a_0 \equiv 2^{2g}\bmod{2^{2g+2}}$ and $a_i\equiv 0\bmod{2^{2g+2-i}}$ for all $i>0$;
 \item $f_0(x)$ has type $1-\{2\}$ at distinct primes $p_t,p_t'>g$;
 \item $f_0(x)$ has $g$ distinct double roots in $\overline{\F}_\ell$ for every odd prime $\ell\leq g$;
 \item $f_0(x)$ has type $1-\{q_1,q_2\}$ at a prime $p_2>2g+2$, which is a primitive root modulo~$q_1, q_2$ and $q_3$, and $p_2\equiv 1 \bmod{3}$;
 \item $f_0(x)$ has type $2{-}\{q_3\}$ at a prime $p_3>2g+2$, which is a primitive root modulo~$q_3$, and $p_3\equiv 1 \bmod{3}$;
 \item $f_0(x)$ has type $1-\{q_4,q_5\}$ at a prime $p_2'>2g+2$, which is a primitive root modulo~$q_3,q_4$ and $q_5$;  
 \item $f_0(x)$ has type $2{-}\{q_5\}$ at a prime $p_3'>2g+2$, which is a primitive root modulo~$q_5$;
 \item $f_0(x)$ modulo a prime $p_{irr}$ is irreducible ;
 \item $f_0(x)$  modulo a prime $p_{lin}$ factors as an irreducible polynomial times a linear factor.
\end{itemize}
Let $C:y^2=f(x)$ be a hyperelliptic curve over $\Q$ with $f(x)\in \Z[x]$ monic and squarefree such that 
\begin{enumerate}
 \item $f(x)\equiv f_0\bmod{N}$, where
  $$N=p_t^2\cdot p_t'^2\cdot p_{lin}\cdot p_{irr}\cdot p_2^2\cdot p_2'^2\cdot p_3^3\cdot p_3'^3\cdot 2^{2g+2}\cdot \prod_{\footnotesize{\shortstack{$3\leq p \leq g$ \\prime}}} p^{2},$$
 \item $f(x)$ mod~$p$ has no roots of multiplicity greater than $2$ in $\overline{\F}_p$ for all primes $p$ not dividing $N$.
\end{enumerate}
Then $$\Gal(\Q(J[\ell])/\Q)\cong \begin{cases} 
\GSp_{2g}(\F_\ell) \mbox{ for all primes }\ell\neq 2,\\
S_{2g+2} \mbox{ for }\ell=2,
\end{cases}$$
where $J=\Jac(C)$.\end{thm}

\proof
Clearly hypothesis $(2\mathrm{G}+\epsilon)$ of Theorem~\ref{alldone} is satisfied. 

Since $f(x)\equiv f_0\bmod{N}$ then by Lemma~\ref{type_cong} $f(x)\in \Z[x]$ satisfies hypotheses $(2\mathrm{T})$, $(p_2)$, $(p_3)$, $(p_2')$ and $(p_3')$ of Theorem~\ref{alldone}.

Hypotheses $(\mathrm{TT})$ and $(\mathrm{ss})$ are satisfied too by Lemma~\ref{roots:multiplicity2}, Corollary~\ref{explicit_tt} and  Lemma~\ref{good:2}~(ii).

Hypothesis $(3)$ holds since $p_2,p_3\equiv 1 \bmod{3}$. 

The existence of $p_{irr}$ and  $p_{lin}$ guarantees that $(S_{2g+2})$ is satisfied.

Therefore by Theorem~\ref{alldone} we have that $\Gal(\Q(J[\ell])/\Q)\cong\GSp_{2g}(\F_\ell)$ for all odd primes and $\Gal(\Q(J[2])/\Q)\cong S_{2g+2}$.
\endproof

\begin{rmk}
\label{explicit_notriple}
Condition $(2)$ can be made explicit, in the sense that one can construct examples for $(2)$ in a systematic way as follows.

Recall that  $f(x) \bmod{p}$ has a root of multiplicity greater than $2$ if and only if $f, f', f'' \bmod{p}$ have a common root in $\overline{\F}_p$.
To construct a suitable polynomial, first pick any $f(x)$ satisfying $(1)$ and such that $f(x) \bmod{p}$ has no roots of multiplicity greater than $2$ for all primes $p<2g$ not dividing $N$. 
Let $$\tilde{N}=N\cdot\prod_{\footnotesize{\shortstack{$p < 2g,\; p\nmid N$ \\prime}}} p.$$
By changing the linear term of $f(x)$ by a multiple of $\tilde{N}$, ensure that $f'(x)$ and $f''(x)$ have no common roots in $\overline{\Q}$, so that for some polynomials $a(x),b(x)\in\Z[x]$ we have $a(x)f'(x)+b(x)f''(x)=M\in \Z\setminus\{0\}$. 
This guarantees that $f(x) \bmod{p}$ does not have roots of multiplicity greater than $2$ for all primes $p \nmid M$. 

If $p\mid M$ with $p\nmid \tilde{N}$ then there exists $c\in \F_p$ such that $f(x)+c  \bmod p$ is non-zero at the $\overline{\F}_p$-roots of $f''(x)$, as $p>2g=\deg f''$. 
Thus, by the Chinese Remainder Theorem, there exist $z\in \Z$ such that $f(x)+z\tilde{N} \bmod{p}$ is non-zero at the $\overline{\F}_p$-roots of $f''(x)$ for every $p\mid M$ with $p\nmid \tilde{N}$. Hence, $f(x)+z\tilde{N}$ satisfies conditions~$(1)$ and~$(2)$ as required.
\end{rmk}

\medskip
We now turn to the proof of the congruence conditions used in the proof of Theorem~\ref{final:congruence_th}. For the remainder of this section $F$ will be a local field of odd residue characteristic~$p$. 
Let $C:y^2=f(x)$ be a hyperelliptic curve over $F$ with $f(x)\in \O_F[x]$ monic and squarefree and let $J=\Jac(C)$. 

\subsection{Congruences and type $t-\{q_1,\dots,q_k\}$}$\;$\\
\vspace{-.3cm}

The description of polynomials of type $t-\{q_1,\dots,q_k\}$ in terms of congruences follows from the following version of Hensel's lemma for lifting factorisations (see  {\cite[ III.4.3, Th\'eor\`eme~1]{Bourbaki_CA}}):


\begin{thm}[Hensel's Lemma for lifting factorisations]
\label{hensel}
Let $F$ be a local field and let $f(x)\in \O_F[x]$ be a monic polynomial. Let $m\geq1$ and suppose that $$f(x)\equiv \prod_{0\leq i\leq k} g_i(x) \bmod{\pi_F^{m}},$$ 
where $g_i(x)\in \O_F[x]$ are monic polynomials such that for every $i\neq j$ the roots of $\overline{g}_i(x)$ are distinct from the roots of $\overline{g}_j(x)$. Then there exist unique monic polynomials  $\tilde{g}_0(x), \dots, \tilde{g}_k(x)\in \O_F[x]$ such that
$\tilde{g}_i(x)\equiv {g}_i(x)\bmod {\pi_F^{m}}$ and 
$$f(x)= \prod_{0\leq i\leq k} \tilde{g}_i(x).$$
\end{thm}

\begin{lem}
\label{type_cong}
Let $f_0(x), f(x)\in \O_F[x]$ be monic polynomials. If $f_0(x)$ has type $t-\{q_1,\dots,q_k\}$ and 
$$f(x)\equiv f_0(x) \bmod{\pi_F^{t+1}},$$ then $f(x)$ has type $t-\{q_1,\dots,q_k\}$.  
\end{lem}

\proof The result follows from Theorem~\ref{hensel} with $m=t+1$ by Definition~\ref{type} and Definition~\ref{t-eis}.
\endproof

\subsection{Semistability at odd primes}

\begin{lem}
\label{roots:multiplicity2}
Suppose $p$ is an odd prime and $f(x)\in \O_F[x]$ is a monic polynomial.
\begin{enumerate}[$(i)$]
 \item If  all roots of $\overline{f}(x)$ in $\overline{\F}_p$ have multiplicity at most $2$, then $J$ is semistable. Moreover, if there are $d$ roots of multiplicity $2$, then $J$ has toric dimension $\min(d,g)$.
 \item If $\overline{f}(x)$ is separable or $f(x)\in \O_F[x]$ has type $t-\{2,2,{\dots},2\}$, where the number of twos is between $1$ and $g+1$, then $J$ is semistable. 
\end{enumerate}
\end{lem}

\proof Clearly $(ii)$ follows from $(i)$.

For simplicity we will use the results and notation of Section~\ref{section_inertia_clusters} to prove $(i)$. 

Let $R$ be the set of roots of $f(x)$, with $\alpha_1, \alpha_1', \dots, \alpha_d, \alpha_d'$ the roots that reduce to double roots in $\overline{\F}_p$, 
i.e.~ $\overline{\alpha}_i=\overline{\alpha}_i'$. The clusters are singleton roots, the set $R$ and $\s_i=\{\alpha_i,\alpha_i'\}$ for $i=1,\dots,d$. 
We readily compute
$$d_R=0,\quad\quad \mu_R=\lambda_R=0, \quad\quad  \epsilon_R=\gamma_R=\triv,  \quad\quad \mathrm{ and\; so\;\; }  V_R=\triv^{\oplus(2g-2d)};$$
and
$$I_{\s_i}=I_F,\quad\quad \mu_{\s_i}=0, \quad\quad \lambda_{\s_i}=d_{\s_i}\in \frac{1}{2}\Z, \quad\quad \epsilon_{\s_i}=\triv,$$  
$$\gamma_{\s_i}=\begin{cases}
 \triv & \text{ if } \lambda_{\s_i}\in \Z,\\
 \text{order two } & \text{ if } \lambda_{\s_i}\not\in\Z,\end{cases} \quad\quad \mathrm{ and\; so\;\; } V_{\s_i}=0.$$
It follows from Theorem~\ref{cluster:theorem} that
$$H^1_{\acute{e}t}(C/\overline{F}, \Q_\ell) = \begin{cases} \triv^{\oplus (2g-2d)} \oplus (\triv^{\oplus d}\otimes \Sp(2)) & \text{ if }d<g+1,\\
\triv^{\oplus g}\otimes \Sp(2) & \text{ if }d=g+1;\end{cases}$$
where $\ell$ is any prime $\ell\neq p$.
In particular, inertia acts unipotently on $H^1_{\acute{e}t}(C/\overline{F}, \Q_\ell)$, so $J$ is semistable (see \cite[Proposition~3.5]{Grothendieck1972}) and has toric dimension $\min(d,g)$.\endproof

\begin{cor}
\label{explicit_tt}
Let $p$ be an odd prime and suppose that $\overline{f}(x)$ has $g$ double roots over $\overline{\F}_p$. Then $J$ is semistable and has totally toric reduction. 
\end{cor}


\subsection{Good reduction at $p=2$}

\begin{lem}
\label{good:2}
Let $F$ be a finite extension of $\Q_2$ and let 
$$f(x)=x^{2g+2}+a_{2g+1}x^{2g+1}+\dots+a_1 x+a_0 \in F[x].$$ 
If either
\begin{enumerate}[(i)]
 \item $a_0-\frac{1}{4}\in \O_F, a_{2g+1} \in \O_F^\times$ and $a_i\in \O_F$ for $1\leq i \leq 2g$, or
 \item $a_0 \equiv 2^{2g}\bmod{2^{2g+2}}, a_{2g+1} \equiv 2\bmod 4$ and $a_i\equiv 0\bmod{2^{2g+2-i}}$ for $1\leq i \leq 2g$,
\end{enumerate}
then $C$ has good reduction. In particular $I_F$ acts trivially on $J[\ell]$ for every odd prime $\ell$.
\end{lem}

\proof (ii) The substitution $x=2X$, $y=2^{g+1}Y$ shows that (i) implies (ii).

(i)  The substitution $y=Y+\frac{1}{2}$ transforms the model of $C$ into 
$$Y^2+Y=f(x)-\frac{1}{4}\in\O_F[x].$$
All points on this affine chart are smooth since the partial derivative with respect to $Y$ is nowhere vanishing. 
The substitution $V=1/x, W=Y/x^{g+1}$ gives the chart at infinity $W^2 + WV^{g+1} = V^{2g+2} (f(1/V)-\frac{1}{4})$. 
There is a unique point at infinity, corresponding to $V=0$, which is a smooth point since the partial derivative with respect to $V$ is a unit: the linear term of the RHS is $a_{2g+1}V$. 
Therefore, the curve has good reduction at $2$.

The last statement then follows from the theorem of N\'eron\--Ogg\--Shafarevich.\endproof

\section{An example}\label{section_example}

In this section we construct an explicit hyperelliptic curve of genus $6$ with maximal mod~$\ell$ Galois representation for all primes $\ell$, following the recipe of Theorem~\ref{final:congruence_th}.
 
First of all, $2g+2=14=7+7=3+11$, so we can take $q_1=q_2=7$, $q_4=3$, $q_5=11$ and $q_3=13$. 
Now pick primes that satisfy the appropriate congruence conditions:
$$p_t=7,\; p_t'=11,\; p_{lin}=23,\; p_{irr}=29, \; p_2=19, \; p_2'=41, \; p_3=37, \; p_3'=17.$$
For example $p_2$ is a primitive root modulo $7$ and $13$ and it is congruent to $1$ modulo $3$, so the choice $p_2=19$ meets the requirements of Theorem~\ref{final:congruence_th}, and similarly for the other primes.

The theorem then gives the following requirements for 
$$f_0(x)=x^{14}+a_{13}x^{13}+\dots+a_1 x+a_0 \in \Z[x]:$$ 
{\footnotesize
\noindent
\makebox[5.5cm][l]{$f_0(x)$ has type $1-\{2\}$ at $7$,}  $f_0(x)$ has type $1-\{2\}$ at $11$, \\
\noindent
\makebox[5.5cm][l]{$f_0(x)$ has type $1-\{7,7\}$ at $19$,}  $f_0(x)$ has type $1-\{3,11\}$ at $41$,\\
\noindent
\makebox[5.5cm][l]{$f_0(x)$ has type $2-\{13\}$ at $37$,} $f_0(x)$ has type $2-\{11\}$ at $17$,\\
\noindent
\makebox[5.5cm][l]{$f_0(x)$ is irreducible mod~$23$,} $f_0(x)$ factors as linear $\cdot$ irreducible mod~${29}$,\\
$f_0(x)$ has $6$ distinct double roots over $\overline{\F}_3$ and $\overline{\F}_5$,\\
\noindent
$a_0 \equiv 2^{12}\bmod{2^{14}}$, $a_{13}\equiv 2 \bmod{4}$ and $a_i\equiv 0\bmod{2^{14-i}}$ for  $1\leq i \leq 12$.\\}

\noindent
By Definition~\ref{type}, Lemma~\ref{type_cong} and Corollary~\ref{explicit_tt}, it is enough to have:
{\footnotesize \begin{align*}
f_0(x) &\equiv (x^{12} + 2x^8 + 5x^7 + 3x^6 + 2x^5 + 4x^4 + 5x^2 + 3)\cdot (x^2-7) && \bmod 7^2,\\
f_0(x) &\equiv (x^{12} + x^{8} + x^{7} + 4 x^{6} + 2 x^{5} + 5 x^{4} + 5 x^{3} + 6 x^{2} + 5 x + 2)\cdot (x^2-11)  && \bmod 11^2,\\
f_0(x) &\equiv (x^7-19)\cdot((x-1)^7-19)  && \bmod 19^2,\\
f_0(x) &\equiv (x^{11}-41)\cdot((x-1)^3-41)  && \bmod 41^2,\\
f_0(x) &\equiv (x^{13}-37^2)\cdot(x+1)  && \bmod 37^3,\\
f_0(x) &\equiv (x^{11}-17^2)\cdot(x^3 + x + 14)  && \bmod 17^3,\\
f_0(x) &\equiv x^{14} + x^8 + 5x^7 + 16x^6 + x^5 + 18x^4 + 19x^3 + x^2 + 22x + 5  && \bmod 23,\\
f_0(x) &\equiv (x+1)\cdot(x^{13} + 7x + 27)  && \bmod 29,\\
f_0(x) &\equiv (x - 1) \cdot x \cdot (x^{6} + 2 x^{4} + x^{2} + 2 x + 2)^{2}  && \bmod 3^2,\\
f_0(x) &\equiv (x - 1) \cdot x \cdot (x^{6} + x^{4} + 4 x^{3} + x^{2} + 2)^{2}  && \bmod 5^2,\\
f_0(x) &\equiv x^{14}+2x^{13}+2^{12}  && \bmod 2^{14}.
\end{align*}}

\noindent
By the Chinese Remainder Theorem on the coefficients we obtain the following polynomial for $f_0(x)$:\\
{\footnotesize 
$$
\arraycolsep=1.5pt\def\arraystretch{1}
\begin{array}{crrlrrlr}
f_0(x)= &x^{14}\;+& 1122976550518058592759939074 &x^{13} &+& 10247323490706358348644352 &x^{12}&+\\
        &        +& 1120184609916242124087443456 &x^{11} &+& 186398290364786000921886720   &x^{10}&+\\
	&        +& 1685990245699349559300014080 &x^{9}  &+& 387529952672653585935499264   &x^{8}&+\\
	&        +& 1422826957983635547417870336 &x^{7}  &+& 585983998625429997308035072   &x^{6}&+\\
	&        +& 607434202225985243206107136  &x^{5}  &+& 1820210247550502007557029888  &x^{4}&+\\
	&        +& 533014336994715937945092096  &x^{3}  &+& 595803405154942945879752704   &x^{2}&+\\
	&        +& 1276845913825955586899050496 &x      &+& 1323672381818030813822668800. &&
\end{array}$$}

\noindent 
The reduction modulo~$p$ of the polynomial $f_0(x)$ has no roots of multiplicity greater than $2$ for any prime $p\notin\{19, 41, 37, 17\}$, so by 
Theorem~\ref{final:congruence_th} the Jacobian $J_0$ of $$C_0:y^2=f_0(x)$$
has 
$$\Gal(\Q(J_0[\ell])/\Q)\cong \begin{cases} 
\GSp_{12}(\F_\ell) \mbox{ for all primes }\ell\neq 2,\\
S_{14} \mbox{ for }\ell=2.
\end{cases}$$
Moreover, setting 
$$N=p_t^2\cdot p_t'^2\cdot p_{lin}\cdot p_{irr}\cdot p_2^2\cdot p_2'^2\cdot p_3^3\cdot p_3'^3\cdot 2^{2g+2}\cdot 3^{2}\cdot 5^2= 2201590757511816436065484800,$$
the same conclusion holds for any curve $C:y^2=f(x)$ with $f(x)\equiv f_0(x)\bmod N$ such that $f(x)\bmod p$ has no roots of multiplicity greater than $2$ for all primes $p\notin\{19, 41, 37, 17\}$.

\bibliographystyle{annotate}
\bibliography{biblio}

\end{document}